\documentclass[a4paper,twoside]{article}

  \usepackage{amsthm}
  \usepackage{amsmath}
  \usepackage{amsfonts}
  \usepackage{amssymb}
  \usepackage{amscd}
  \usepackage[mathscr]{eucal}

  \usepackage{enumitem}

  \usepackage{fancyhdr}
  \usepackage[linkcolor=black,colorlinks=true]{hyperref}

  \usepackage[top=2.5cm, bottom=2.7cm,left=3.5cm, right=2.5cm]{geometry}

  \swapnumbers 

  \theoremstyle{definition}  
   \newtheorem{defn}{Definition}[subsection]
   \newtheorem{eg}[defn]{Example}
   \newtheorem*{prblm}{Problem}
   \newtheorem*{ackn}{Acknowledgment}

  \theoremstyle{plain}  
   \newtheorem{thm}[defn]{Theorem}
   \newtheorem{lem}[defn]{Lemma}
   \newtheorem{prop}[defn]{Proposition}
   \newtheorem{cor}[defn]{Corollary}
   \newtheorem*{conj}{Conjecture}

  \theoremstyle{remark} 
   \newtheorem*{note}{Note}
   \newtheorem{rmk}[defn]{Remark}
   \newtheorem*{notn}{Notation}

   \newtheorem*{claim}{Claim}

  \numberwithin{equation}{section}
  \allowdisplaybreaks[1] 
  \setlist[enumerate]{font=\upshape} 

  \newcommand{\bdefn}{\begin{defn}}
  \newcommand{\beg}{\begin{eg}}
  \newcommand{\bprblm}{\begin{prblm}}
  \newcommand{\blem}{\begin{lem}}
  \newcommand{\bprop}{\begin{prop}}
  \newcommand{\bthm}{\begin{thm}}
  \newcommand{\bcor}{\begin{cor}}
  \newcommand{\bconj}{\begin{conj}}
  \newcommand{\brmk}{\begin{rmk}}
  \newcommand{\bnote}{\begin{note}}
  \newcommand{\bnotn}{\begin{notn}}
  \newcommand{\bclaim}{\begin{claim}}
  \newcommand{\bproof}{\begin{proof}}
  \newcommand{\backn}{\begin{ackn}}

  \newcommand{\edefn}{\end{defn}}
  \newcommand{\eeg}{\end{eg}}
  \newcommand{\eprblm}{\end{prblm}}
  \newcommand{\elem}{\end{lem}}
  \newcommand{\eprop}{\end{prop}}
  \newcommand{\ethm}{\end{thm}}
  \newcommand{\ecor}{\end{cor}}
  \newcommand{\econj}{\end{conj}}
  \newcommand{\ermk}{\end{rmk}}
  \newcommand{\enote}{\end{note}}
  \newcommand{\enotn}{\end{notn}}
  \newcommand{\eclaim}{\end{claim}}
  \newcommand{\eproof}{\end{proof}}
  \newcommand{\eackn}{\end{ackn}}

  \renewcommand{\bf}[1]{\textbf{#1}}
  \renewcommand{\it}[1]{\textit{#1}}
  \renewcommand{\sc}[1]{\textsc{#1}}
  \renewcommand{\sf}[1]{\textsf{#1}}

  \newcommand{\mbb}[1]{\mathbb{#1}}
  \newcommand{\mcl}[1]{\mathcal{#1}}

  \newcommand{\msc}[1]{\mathscr{#1}}
  \newcommand{\msf}[1]{\mathsf{#1}}
  
  \renewcommand{\mit}[1]{\mathit{#1}}

  \newcommand{\ol}[1]{\overline{#1}}

  \newcommand{\wtilde}[1]{\widetilde{#1}}

  \newcommand{\vphi}{\varphi}

  \newcommand{\abs}[1]{\left\lvert#1\right\rvert}
  \newcommand{\babs}[1]{\bigl\lvert#1\bigr\rvert}

  \newcommand{\norm}[1]{\left\lVert#1\right\rVert}
  \newcommand{\bnorm}[1]{\bigl\lVert#1\bigr\rVert}
  
  \newcommand{\snorm}[1]{\norm{\smash{#1}}}

  \newcommand{\B}[1]{\msc{B}({#1})}
  \newcommand{\K}[1]{\msc{K}({#1})}
  
  \newcommand{\BA}[1]{\msc{B}^a({#1})}
  \newcommand{\BL}[1]{\msc{B}^{a,bil}({#1})}

  \newcommand{\ip}[1]{\langle#1\rangle}

  \newcommand{\ran}[1]{\sf{ran}(#1)}
  
  \renewcommand{\Re}[1]{\textnormal{Re}(#1)}

  \DeclareMathOperator{\I}{\sf{1}}

  \DeclareMathOperator{\lspan}{\sf{span}}
  \DeclareMathOperator{\cspan}{\ol{\lspan}}
  \DeclareMathOperator{\sspan}{\ol{\lspan}^{\,\mit{s}}}

  \DeclareMathOperator*{\soplus}{\ol{\oplus}^{\,\mit{s}}}

\begin{document}
  \fontsize{12}{14}
  \selectfont
  \title{\LARGE Bures Distance For Completely Positive Maps}
  \author{B.V. Rajarama Bhat and K. Sumesh}
  \date{}
  \maketitle

 \begin{abstract}
  D. Bures had defined a metric on the set of normal states on a von Neumann algebra using GNS representations of states. This notion has been extended to completely positive maps between $C^*$-algebras by D. Kretschmann, D. Schlingemann and R. F. Werner. We present a Hilbert $C\sp{\ast}$-module version of this theory. We show that we do get a metric when the completely positive maps under consideration map to  a von Neumann algebra. Further, we include several examples and counter examples.
  We also prove a rigidity theorem, showing that representation modules of completely positive maps which are close to the identity map contain a copy of the original algebra.\\

 \end{abstract}


 \section{Preliminaries}
 \subsection{Introduction}
  Given a state $\phi $ on a $C^*$-algebra $\mcl{B}$ we have the familiar GNS-triple $(H,\pi,x)$, where $H$ is a Hilbert space, $\pi:\mcl{B}\rightarrow\B{H}$ is a representation (i.e., $\ast$-homomorphism)  and $x\in H$ is a vector such that $\phi(\cdot)=\ip{x,\pi(\cdot)x}$. Now it is a natural question to ask: If two states $\phi_1,\phi_2$ are close in some metric, whether the associated triples are close in some sense? Keeping this idea in mind, D. Bures (\cite{Bures}) defines a distance between two states $\phi_1,\phi_2$ on $\mcl{B}$, as
  \begin{equation*}
    \beta(\phi_1,\phi_2):=\inf\norm{x_1-x_2},
  \end{equation*}
  where the infimum is taken over all GNS-triples with common representation spaces: $(H,\pi,x_1), (H,\pi,x_2)$ of $\phi_1,\phi_2$. D. Bures showed that $\beta$ is indeed a metric. The notion has
  found uses in many areas  (\cite{Ab, alberti, araki, Ditt98}).

  D. Kretschmann, D. Schlingemann and R. F. Werner (\cite{KDR08}) extended this notion at first to completely positive (CP) maps from a $C^*$-algebra $\mcl{A}$ to $\B{G}$ for some Hilbert space $G$ and then to more general range $C^*$-algebras using an alternative definition of the Bures distance. They use the Stinespring representation (\cite{wfs}) for the initial definition, which in the usual formulation requires the range space to be the whole algebra $\B{G}$. Here we develop the theory using Hilbert module language, which allows the range algebra to be any $C^*$-algebra, and the definition of the metric is a very natural extension of the definition given by Bures for states. Working with modules has several advantages.  The results we get are of course same as that of \cite{KDR08}, when the range algebra is a von Neumann algebra or an injective $C^*$-algebra. However,  we show that one may not even get a metric (triangle inequality may fail) when the range algebra is a general $C^*$-algebra.

  There have been several papers (\cite{Akh, Ditt99}) on different methods to make exact computations of the Bures metric for states. We provide several examples with explicit computations of the Bures distance for CP-maps. In
  particular, we show that the infimum in the definition of Bures metric may not be attained in all common representation modules, answering a question raised in \cite{KDR06,KDR08}.
  It turns out that the example is quite simple involving CP-maps on $2\times 2$ matrix algebra.

  In the last Section we prove a rigidity theorem, which says that on von Neumann algebras, if a CP-map is strictly within unit distance (in Bures metric) from  the identity map, then the GNS-module of the CP-map contains a copy of the original von Neumann algebra as a direct summand. We consider this as the most important positive result of this paper and we expect that the result will have further applications in the study of CP-maps.

 \subsection{Hilbert \texorpdfstring{$C^*$}{C*}-modules}
  Let $\mcl{B}$ be a $C^*$-algebra. A complex vector space $E$ is a Hilbert $\mcl{B}$-module if it is a right $\mcl{B}$-module with a $\mcl{B}$-valued inner product, which is complete with respect to the associated norm (see \cite{lance, paschke73, skeide00} for basic theory).  We denote the space of all bounded and adjointable maps between two Hilbert $\mcl{B}$-modules $E_1$ and $E_2$ by $\BA{E_1,E_2}$.  In particular, if $E_1=E_2=E$, then $\BA{E,E}=\BA{E}$, which forms a $C^*$-algebra with natural algebraic operations.

  Let  $\pi:\mcl{B}\rightarrow\B{G}$ be a non-degenerate (i.e., $\cspan\pi(\mcl{B})G=G$) representation of $\mcl{B}$ on a Hilbert space $G$. Given a Hilbert $\mcl{B}$-module $E$,  we define the Hilbert space $H:=E\odot G$ as the completion of the inner product space obtained from the algebraic tensor product $E\otimes G$ by dividing out the null space of the semi-inner product $\ip{x\otimes g,x'\otimes g'}:=\ip{g,\pi(\ip{x,x'})g'}$, where $x,x'\in E,~g,g'\in G$. We denote the equivalence class containing $x\otimes g$ by $x\odot g$. To each $x\in E$ we associate the linear map $L_x:g\mapsto x\odot g$ in $\B{G,H}$ with adjoint $L_x^*:y\odot g\mapsto \pi(\ip{x,y})g$. Clearly  $L_x^*L_y=\pi(\ip{x,y})$ and $L_{xb}=L_x\pi(b)$ for all $x,y\in E,~b\in \mcl{B}$. Also $\norm{L_x}^2=\norm{\pi(\ip{x,x})}=\norm{x}^2$. By identifying  $\mcl{B}$ with $\pi(\mcl{B})$ and $x$ with $L_x$ we may assume that  $E\subseteq\B{G,H}$. Note that $a\mapsto a\odot id_G:\BA{E}\rightarrow \B{H}$ is a unital $\ast$-homomorphism, and hence  an isometry. So we may consider $\BA{E}\subseteq \B{H}$.

  Suppose $\mcl{A}$ is another $C^*$-algebra. A Hilbert $\mcl{B}$-module $E$ is called a Hilbert $\mcl{A}$-$\mcl{B}$-module if there exists a  representation $\tau:\mcl{A}\rightarrow \BA{E}$ which is non-degenerate (equivalently,
  unital if $\mcl{A}$ is unital). If $E$ is a Hilbert $\mcl{A}$-$\mcl{B}$-module, we may consider $\mcl{A}\subseteq \BA{E}$, and we denote $\tau(a)$ by $a$ itself and thereby $\tau(a)x=ax$ for all $x\in E,~a\in \mcl{A}$. The composition map $\mcl{A}\xrightarrow{\tau} \BA{E}\rightarrow \B{H}$ is denoted by $\rho$; i.e., $\rho(a)=a\odot id_G$. Note that $L_{ax}=\rho(a)L_x$.  Also $\B{G,H}$ forms a Hilbert $\mcl{A}$-$\B{G}$-module with left action $ax:=\rho(a)x$. If $E_1$ and $E_2$ are two Hilbert $\mcl{A}$-$\mcl{B}$-modules,  then a linear map $\Phi:E_1\rightarrow E_2$ is said to be $\mcl{A}$-$\mcl{B}$-linear (or bilinear) if $\Phi(axb)=a\Phi(x)b$ for all $a\in\mcl{A},~b\in\mcl{B},~x\in E$. The space of all bounded, adjointable and bilinear maps from $E_1$ to $E_2$ is denoted by $\BL{E_1,E_2}$. If $E$ is a Hilbert $\mcl{A}$-$\mcl{B}$-module,  then $\BL{E}$ is the relative commutant of the image of $\mcl{A}$ in $\BA{E}$.

  Suppose $\mcl{B}\subseteq \B{G}$ is a von Neumann algebra and $E$ is a Hilbert $\mcl{B}$-module. Then we say $E$ is a  von Neumann $\mcl{B}$-module if $E$ is strongly closed in $\B{G,H}\subseteq \B{G\oplus H}$. Thus, if $x$ is an element in the strong closure $\ol{E}^s$ of a  Hilbert $\mcl{B}$-module  $E$, then there exists a net $(x_{\alpha})\subseteq E$ such that $L_{x_{\alpha}}\xrightarrow{SOT} x$. All von Neumann $\mcl{B}$-modules are self-dual (in the sense that all $\mcl{B}$-valued functionals are given by a $\mcl{B}$-valued inner product), and hence they are complemented in all Hilbert $\mcl{B}$-module which contains it as a $\mcl{B}$-submodule. In particular, strongly closed $\mcl{B}$-submodules are complemented in a von Neumann $\mcl{B}$-module. If we think $\BA{E}\subseteq \B{H}$, then $\BA{E}$ is a von Neumann algebra acting non-degenerately on the Hilbert space $H$. If $\mcl{A}$ is a $C^*$-algebra, then by a von Neumann $\mcl{A}$-$\mcl{B}$-module we mean a von Neumann $\mcl{B}$-module $E$ with a non-degenerate representation $\tau:\mcl{A}\rightarrow \BA{E}$. In addition, if $\mcl{A}$ is a von Neumann algebra and $a\mapsto \ip{x,ax}:\mcl{A}\rightarrow \mcl{B}$ is a normal mapping for all $x\in E$  (equivalently,  the representation $\rho:\mcl{A}\rightarrow \B{H}$ is normal), then we call $E$ a two-sided von Neumann $\mcl{A}$-$\mcl{B}$-module. For more details see \cite{paschke73,skeide00,skeide05}.

  It is well-known that if $\vphi:\mcl{A}\rightarrow\mcl{B}$ is a CP-map between unital $C^*$-algebras, then there exists a Hilbert $\mcl{A}$-$\mcl{B}$-module $E$ and $x\in E$ such that $\vphi(a)=\ip{x,ax}$ for all $a\in \mcl{A}$. The construction of $E$ is by starting with $\mcl{A}\otimes \mcl{B}$ and defining a $\mcl{B}$-valued semi-inner product on it as $\ip{a_1\otimes b_1,a_2\otimes b_2}:=b_1^*\vphi(a_1^*a_2)b_2$, and usual quotienting and completion procedure (see \cite{kasparov80, Murphy, paschke73, paulsen, wfs}).  The comparison with GNS for states is obvious. The pair $(E,x)$ is called a GNS-construction for $\vphi$ and $E$ is called a GNS-module for $\vphi$. If further, $\cspan\,\mcl{A}x\mcl{B}=E$, then $(E,x)$ is said to be a minimal GNS-construction, and is unique up to isomorphism. If (both $\mcl{A}$ and) $\mcl{B}$ is a von Neumann algebra, then $(E,x)$ can be chosen such that $E$ is a (two-sided) von Neumann $\mcl{A}$-$\mcl{B}$-module. Here the closure for minimality is taken under strong operator topology.

  Note that if $\mcl{B}=\B{G}$, then $L_x^*\rho(a)L_x=\ip{x,ax}=\vphi(a)$ for all $a\in \mcl{A}$. Thus $(\rho, L_x, H)$ is a Stinespring representation for the CP-map $\vphi:\mcl{A}\rightarrow \B{G}$.

 \subsection{Bures distance}
  Given two unital $C^*$-algebras $\mcl{A}$ and $\mcl{B}$, we let $CP(\mcl{A},\mcl{B})$ denote the set of all nonzero CP-maps from $\mcl{A}$ into $\mcl{B}$.

 \bdefn
  A  Hilbert $\mcl{A}$-$\mcl{B}$-module $E$ is said to be a \it{common representation module} for $\vphi_1,\vphi_2\in CP(\mcl{A},\mcl{B})$ if both of them can be represented in $E$, that is, there exist $x_i\in E$ such that $\vphi_i(a)=\ip{x_i,ax_i},~i=1,2$.
 \edefn

  Note that we are demanding no minimality for the common representation module. So we can always have such a module. For, if $(\hat{E_i},\hat{x_i})$ is the minimal GNS-construction for $\vphi_i$, then take $E=\hat{E_1}\oplus\hat{E_2},x_1=\hat{x_1}\oplus 0$ and $x_2=0\oplus\hat{x_2}$. For a common representation module $E$, define $S(E,\varphi_i)$ to be the set of all $x\in E$ such that $\vphi_i(a)=\ip{x,ax}$ for all $a\in\mcl{A}$.

 \bdefn
  Let $E$ be a common representation module for $\vphi_1,\vphi_2\in CP(\mcl{A},\mcl{B})$. Define
  \begin{equation*}
     \beta_E(\vphi_1,\vphi_2):=\inf\big\{\norm{x_1-x_2} :  x_i\in S(E,\vphi_i), i=1,2\big\}
  \end{equation*}
  and the \it{Bures distance}
  \begin{equation*}
     \beta(\vphi_1,\vphi_2):=\underset{E}{\inf}\,\beta_E(\vphi_1,\vphi_2)
  \end{equation*}
  where the infimum is taken over all common representation module $E$.
 \edefn

  We have called $\beta $ as a `distance' in anticipation. Later we will show that it is indeed a metric under most situations, for instance, when $\mcl{B}$ is a von Neumann algebra. But surprisingly $\beta $ is not a metric in general.

  Our first job is to show that the definition here matches with that of \cite{KDR08}.  We see it as follows. Suppose $\mcl{B}=\B{G}$. If $E$ is a common representation module and $x_i\in S(E,\vphi_i)$, then $(\rho,L_{x_i},E\odot G)$ is a Stinespring representation for $\vphi_i$ with $\norm{x_1-x_2}=\norm{L_{x_1}-L_{x_2}}$.  On the other way if
  $(\pi',V_i,H')$ is a  Stinespring representation for $\vphi_i$, then $E:=\B{G,H'}$ is a   Hilbert\footnote{If $\mcl{A}$ is a von Neumann algebra and  $\pi$ is normal, then $a\mapsto\ip{x,ay}=\ip{x,\pi(a)y}$ is normal map from $\mcl{A}\rightarrow\mcl{B}$ for all $x,y\in E$. Thus $E$ can be made into a two-sided von Neumann $\mcl{A}$-$\B{G}$-module.} $\mcl{A}$-$\B{G}$-module with inner product $\ip{x_1,x_2}:=x_1^*x_2$, composition as the right module action and left action given by $ax:=\pi'(a)x$ for all $a\in\mcl{A},~x\in E$. Clearly $(E,V_i)$ is a GNS-construction for $\vphi_i$. Note that $\cspan EG=H'$. We have $H:=E\odot G$ is a Hilbert space with inner product $\ip{x\odot g,x'\odot g'}=\ip{g',x^*x'g'}=\ip{xg,x'g'}$. Thus $x\odot g\mapsto xg$ is a unitary from $U: H\rightarrow H'$. Note that $UL_{V_i}=V_i$ and $U\rho(a)U^*=\pi'(a)$ for all $a\in \mcl{A}$.
  Identifying $H$ with $H'$ through $U$, we get $\pi'=\rho$ and $L_{V_i}=V_i$. Therefore $(\pi,V_i,H')=(\rho,L_{V_i},H)$. Thus there exists a  one-one correspondence between the GNS-constructions $\{(E,x_1), (E, x_2)\}$ and  the Stinespring representations $\{(\pi',V_1,H'), (\pi ', V_2, H')\}$ such that $\norm{x_1-x_2}=\norm{V_1-V_2}$. Hence $\beta(\vphi_1,\vphi_2)$ coincides with the definition given in \cite{KDR08}. In particular, if $\mcl{B}=\B{\mbb{C}}=\mbb{C}$, then $\beta(\vphi_1,\vphi_2)$ is the Bures distance given in \cite{Bures}.

  The following proposition says that $\beta(\varphi_1,\varphi_2)$ coincide with the alternative definition, given in \cite{KDR08}, of Bures distance for CP-maps between arbitrary $C\sp{\ast}$-algebras. Subsequently will not be needing this definition and we present it here for the sake of completeness.

 \bprop
  With notation as above,
  \begin{equation*}
     \beta(\vphi_1,\vphi_2)=\underset{\vphi}{\inf}\,\norm{\vphi_{11}(1)+\vphi_{22}(1)-\vphi_{12}(1)-\vphi_{21}(1)}^{\frac{1}{2}}
  \end{equation*}
  where the infimum is taken over all CP-extensions $\vphi:\mcl{A}\rightarrow M_2(\mcl{B})$ of the form
  $\vphi=\begin{bmatrix}
                         \vphi_{11} & \vphi_{12} \\
                         \vphi_{21} & \vphi_{22}
         \end{bmatrix}$
  with completely bounded maps $\vphi_{ij}:\mcl{A}\rightarrow \mcl{B}$ satisfying $\vphi_{ii}=\vphi_i$.
 \eprop

 \bproof
  Let $E$ be a common representation module and $x_i\in S(E,\vphi_i)$. Define $\vphi:\mcl{A}\rightarrow M_2(\mcl{B})$ by $a\mapsto [\vphi_{ij}(a)]$, where $\vphi_{ij}(a):=\ip{x_i,ax_j}$. Then $\vphi$ is a CP-map with
  \begin{align*}
      \norm{x_1-x_2}^2 &=\norm{\ip{x_1,x_1}+\ip{x_2,x_2}-\ip{x_1,x_2}-\ip{x_2,x_1}} \\
                       &=\norm{\vphi_{11}(1)+\vphi_{22}(1)-\vphi_{12}(1)-\vphi_{21}(1)}.
  \end{align*}
  Since $E$ is arbitrary $\beta(\vphi_1,\vphi_2)\geq\underset{\vphi}{\inf}\,\norm{\vphi_{11}(1)+\vphi_{22}(1)-\vphi_{12}(1)-\vphi_{21}(1)}^{\frac{1}{2}}$. To get the reverse inequality, assume that $\vphi=[\vphi_{ij}]:\mcl{A}\rightarrow M_2(\mcl{B})$ is a CP-map with $\vphi_{ii}=\vphi_i$. Let $(\hat{E},\hat{x})$ be a GNS-construction of $\vphi$. Note that $\hat{E}$ is a Hilbert $\mcl{A}$-$M_2(\mcl{B})$-module. Given $b\in\mcl{B}, x\in \hat{E}$ define $xb:=x(bI)$, where $I\in M_2(\mcl{B})$ is the identity matrix. Under this action $\hat{E}$ becomes a right $\mcl{B}$-module. Now for $x_1,x_2\in \hat{E}$ define $\ip{x_1,x_2}':=\underset{i,j}{\sum}\ip{x_1,x_2}_{ij}$, where $\ip{x_1,x_2}_{ij}$ is the $(i,j)^{\textnormal{th}}$ entry of $\ip{x_1,x_2}\in M_2(\mcl{B})$. Then $\ip{\cdot,\cdot}'$ is a $\mcl{B}$-valued inner product on $\hat{E}$. Denote  the resulting inner product $\mcl{B}$-module by $E_0$. The left action of $\mcl{A}$ on $\hat{E}$ induce a non-degenerate left action of $\mcl{A}$ on $E_0$. Complete $E_0$ to get the Hilbert $\mcl{A}$-$\mcl{B}$-module $E$. Set $x_i=\hat{x}e_{ii}$,  where
  $\{e_{ij}\},~1\leq i,j\leq 2$ are matrix units of $M_2(\mcl{B})$. Then $x_i\in S(E,\vphi_i)$ and
  \begin{align*}
      \norm{x_1-x_2}^2 &=\norm{\ip{x_1-x_2,x_1-x_2}'} \\
                       &=\norm{\ip{x_1,x_1}'+\ip{x_2,x_2}'-\ip{x_1,x_2}'-\ip{x_2,x_1}'}  \\
                       &=\norm{\vphi_{11}(1)+\vphi_{22}(1)-\vphi_{12}(1)-\vphi_{21}(1)}.
  \end{align*}
  Since $\vphi$ is arbitrary $\beta(\vphi_1,\vphi_2)\leq\underset{\vphi}{\inf}\,\norm{\vphi_{11}(1)+\vphi_{22}(1)-\vphi_{12}(1)-\vphi_{21}(1)}^{\frac{1}{2}}$.
 \eproof

  The following proposition says that Bures distance is stable under taking ampliations.

 \bprop\label{betaampl}
  Let $\mcl{A}$ and $\mcl{B}$ be unital $C^*$-algebras. Then for $\vphi,\psi\in CP(\mcl{A},\mcl{B})$,
  \begin{equation*}
       \beta(\vphi,\psi)=\beta(\vphi_n,\psi_n)
  \end{equation*}
  where $\vphi_n,\psi_n:M_n(\mcl{A})\rightarrow M_n(\mcl{B})$ are the amplifications of $\vphi,\psi$ respectively for $n\geq 1.$
 \eprop

 \bproof
  Fix $n\geq 1.$  Suppose $E$ is a common representation module for $\vphi, \psi$ and $x_1\in S(E,\vphi),x_2\in S(E,\psi)$. Then
   $\textrm{diag}(x_1,\cdots,x_1)\in S(M_n(E),\vphi_n)$ and $\textrm{diag}(x_2,\cdots,x_2)\in S(M_n(E),\psi_n)$,
  and hence
  \begin{equation*}
     \beta(\vphi_n,\psi_n)\leq\norm{\textrm{diag}(x_1-x_2,\cdots,x_1-x_2)}=\norm{x_1-x_2}.
  \end{equation*}
  Since $x_1, x_2$ and $E$ are arbitrary $\beta(\vphi_n,\psi_n)\leq \beta(\vphi,\psi)$. Conversely, suppose $F$ is a common representation module for $\vphi_n,\psi_n$ and $y_1 \in S(F,\varphi_n), y_2\in S(F,\psi_n).$ If $\{e_{ij}\}, \{f_{ij}\},~1\leq i,j \leq n$ are matrix units of $M_n(\mcl{A}), M_n(\mcl{B})$ respectively, then $E:=\{e_{11}Ff_{11}\}$ is a common representation module for $\vphi,\psi$ in the natural way and more over, $e_{11}y_1f_{11}\in S(E,\varphi)$ and $e_{11}y_2f_{11}\in S(E,\psi)$. Also,
  \begin{align*}
     \norm{e_{11}y_1f_{11}-e_{11}y_2f_{11}}^2
               = \norm{f_{11}\ip{e_{11}(y_1-y_2), e_{11}(y_1-y_2)}f_{11}}
               \leq \norm{y_1-y_2}^2.
  \end{align*}
  Therefore $\beta (\vphi,\psi)\leq\beta(\vphi_n,\psi_n).$
 \eproof

 \bprop\label{betacomp}
  Let $\mcl{A},\mcl{B}$ and $\mcl{C}$ be unital $C^*$-algebras. Then for $\vphi_i\in CP(\mcl{A},\mcl{B})$ and $\psi_i\in CP(\mcl{B},\mcl{C}),~i=1,2$,
  \begin{equation*}
        \beta(\psi_1\circ\vphi_1,\psi_2\circ\vphi_2)
                 \leq \norm{\vphi_1}^{\frac{1}{2}}\beta(\psi_1,\psi_2)+ \norm{\psi_2}^{\frac{1}{2}}\beta(\vphi_1,\vphi_2).
  \end{equation*}
  In particular,
  \begin{equation*}
        \beta(\psi_2\circ\vphi_1,\psi_2\circ\vphi_2)
                 \leq  \norm{\psi_2}^{\frac{1}{2}}\beta(\vphi_1,\vphi_2).
  \end{equation*}
 \eprop

 \bproof
  Suppose $E,F$ are common representation modules for $\vphi_i,\psi_i$ respectively, and $x_i\in S(E,\vphi_i),y_i\in S(F,\psi_i),~i=1,2$. Then $x_i\odot y_i\in S(E\odot F,\psi_i\circ\vphi_i)$, and hence
  \begin{align*}
       \beta(\psi_1\circ\vphi_1,\psi_2\circ\vphi_2)
                    & \leq \norm{x_1\odot y_1-x_2\odot y_2}  \\
                    & \leq \norm{x_1\odot y_1-x_1\odot y_2+x_1\odot y_2-x_2\odot y_2}  \\
                    & \leq \norm{x_1}\norm{y_1-y_2}+\norm{x_1-x_2}\norm{y_2}  \\
                    & \leq \norm{\vphi_1}^{\frac{1}{2}} \norm{y_1-y_2}+\norm{x_1-x_2}\norm{\psi_2}^{\frac{1}{2}}.
  \end{align*}
  Since $x_i,y_i,E$ and $F$ are arbitrary the results holds.
 \eproof

 \bprop
  Let $\vphi_1,\vphi_2\in CP(\mcl{A},\mcl{B})$. Then
  \begin{enumerate}
    \item [(i)] $\beta(\vphi_1,\vphi_1+\vphi_2)\leq \norm{\vphi_2}^{\frac{1}{2}}$.
    \item [(ii)] $\babs{\beta(\vphi_1,\vphi_2)-\beta(\vphi_1,\epsilon\vphi_1+(1-\epsilon)\vphi_2)}\leq \epsilon^{\frac{1}{2}}(\norm{\vphi_1}^{\frac{1}{2}}+\norm{\vphi_2}^{\frac{1}{2}})$ for $0\leq\epsilon\leq 1$.
    \item [(iii)] If $\vphi_i(1)\leq 1$, then $\babs{\norm{\vphi_1}-\norm{\vphi_2}}\leq 2\beta(\vphi_1,\vphi_2)$.
  \end{enumerate}
 \eprop

 \bproof
  (i) Suppose $(E_i,x_i)$ is a GNS-construction for $\vphi_i,~i=1,2$. Then $ z_1:=x_1\oplus 0\in S(E_1\oplus E_2,\vphi_1)$ and $z_2:=x_1\oplus x_2\in S(E_1\oplus E_2,\vphi_1+\vphi_2)$, and hence
  \begin{equation*}
     \beta(\vphi_1,\vphi_1+\vphi_2) \leq \norm{z_1-z_2}
                                    =\norm{x_2}
                                    =\norm{\vphi_2}^{\frac{1}{2}}.
  \end{equation*}

  (ii) Using triangle inequality and part (i),
  \begin{align*}
      \babs{\beta(\vphi_1,\vphi_2)-\beta(\vphi_1,\epsilon\vphi_1&+(1-\epsilon)\vphi_2)} \\
                  &\leq\beta(\vphi_2,\epsilon\vphi_1+(1-\epsilon)\vphi_2) \\
                  &\leq\beta(\vphi_2,(1-\epsilon)\vphi_2)+\beta((1-\epsilon)\vphi_2,\epsilon\vphi_1+(1-\epsilon)\vphi_2) \\
                  &\leq\norm{\epsilon\vphi_2}^{\frac{1}{2}}+\norm{\epsilon\vphi_1}^{\frac{1}{2}} \\
                  &\leq \epsilon^{\frac{1}{2}}(\norm{\vphi_1}^{\frac{1}{2}}+\norm{\vphi_2}^{\frac{1}{2}}).
  \end{align*}

  (iii) Let $E$ be a common representation for $\vphi_1,\vphi_2$ and $x_i\in S(E,\vphi_i)$. Then
  \begin{align*}
      \babs{\norm{\vphi_1}-\norm{\vphi_2}} 
                                            &= \babs{\norm{x_1}^2-\norm{x_2}^2} \\
                                            &= \babs{(\norm{x_1}+\norm{x_2})(\norm{x_1}-\norm{x_2})} \\
                                            &= (\norm{x_1}+\norm{x_2})\babs{\norm{x_1}-\norm{x_2}}  \\
                                            &\leq 2 \norm{x_1-x_2}.
  \end{align*}
  Since $x_1,x_2$ and $E$ are arbitrary the result follows.
 \eproof

 \section{Bures distance: von Neumann algebras}
  As is well-known one of the problems in dealing with Hilbert $C^*$-modules in contrast to Hilbert spaces is that in general submodules are not complemented, that is, there is a problem in taking orthogonal complements and writing the whole space as a direct sum. This problem is not there for von Neumann modules.
  Here we generalize almost all the results of \cite{KDR08}, where the results stated mainly for the case when the range algebra is the algebra of all bounded operators on a Hilbert space. The proofs are similar, though we have also taken some ideas from \cite{Bures}. We also give several examples and answer a question of \cite{KDR08} in the negative.

  In this Section we assume that $\mcl{A}$ is a unital $C^*$-algebra, $\mcl{B}\subseteq \B{G}$ is a von Neumann algebra and $\vphi_1,\vphi_2\in CP(\mcl{A},\mcl{B})$.

 \subsection{Metric property}
  To begin with we have the following proposition.

  \bprop
  If $\mcl{B}\subseteq \B{G}$ is a von Neumann algebra, then
  \begin{equation}\label{newbeta}
      \beta(\vphi_1,\vphi_2)=\underset{\msf{E}}{\inf}\,\beta_{\msf{E}}(\vphi_1,\vphi_2)
  \end{equation}
  where the infimum is taken over all common representation modules $\msf{E}$ which are von Neumann $\mcl{A}$-$\mcl{B}$-module.
 \eprop

 \bproof
  Since von Neumann $\mcl{B}$-modules are Hilbert $\mcl{B}$-modules we have $\beta(\vphi_1,\vphi_2)\leq\inf\,\beta_{\msf{E}}(\vphi_1,\vphi_2)$. To get the reverse inequality, assume that $E$ is a  common representation module for $\vphi_1,\vphi_2$. Then $\msf{E}:=\ol{E}^s\subseteq \B{G,E\odot G}$ forms a von Neumann $\mcl{A}$-$\mcl{B}$-module. Since $E\subseteq \msf{E}$ we have $\msf{E}$ is a common representation module for $\vphi_1,\vphi_2$, and hence $\inf\,\beta_{\msf{E}}(\vphi_1,\vphi_2)\leq \beta(\vphi_1,\vphi_2)$.
 \eproof


  As we have taken $\mcl{B}$ as von Neumann algebra for this Section, we may use (\ref{newbeta}) as the definition of Bures distance. Also by a common representation module and GNS-module we will mean a von Neumann $\mcl{A}$-$\mcl{B}$-module. However, note that for all the results here, the algebra $\mcl{A}$ can be a general $C^*$-algebra and the left action by $\mcl{A}$ need not be normal. So we do not need that $\vphi_1,\vphi_2$ to be normal.

  The following result shows the existence of a sort of universal module where we can take infimum to compute the Bures distance.

 \bprop\label{universalbeta}
  There exists a von Neumann $\mcl{A}$-$\mcl{B}$-module $\mcl{E}$ such
  that:
  \begin{enumerate}
    \item [(i)] For all $\vphi_1,\vphi_2\in CP(\mcl{A},\mcl{B})$,  $\beta(\vphi_1,\vphi_2)=\beta_{\mcl{E}}(\vphi_1,\vphi_2)$;
    \item [(ii)] For a fixed $\vphi_1\in CP(\mcl{A},\mcl{B})$ there exists $\xi_1\in S(\mcl{E},\vphi_1)$ such that $\beta(\vphi_1,\vphi_2)=\inf\big\{\norm{\xi_1-\xi_2} : \xi_2\in S(\mcl{E},\vphi_2)\big\}$ for all $\vphi_2\in CP(\mcl{A},\mcl{B})$.
    %
  \end{enumerate}
 \eprop

 \bproof
  For each $\vphi\in CP(\mcl{A},\mcl{B})$ fix a GNS-construction $(E_{\vphi},x_{\vphi})$. Set $H_{\vphi}=E_{\vphi}\odot G$ and $H=\oplus H_{\vphi}$. Then $\mcl{E}_0:=\soplus E_{\varphi}\subseteq \B{G,H}$ is a von Neumann $\mcl{A}$-$\mcl{B}$-module. Note that $S(\mcl{E}_0,\vphi)$ is nonempty for all $\vphi\in CP(\mcl{A},\mcl{B})$. Take $\mcl{E}=\mcl {E}_0\oplus\mcl{E}_0$ which is a von Neumann $\mcl{A}$-$\mcl{B}$-module.

  (i) Suppose $\vphi_1,\vphi_2\in CP(\mcl{A},\mcl{B})$ and $E$ is a common representation module. We will prove that $\beta_{\mcl{E}}(\vphi_1,\vphi_2)\leq\beta_E(\vphi_1,\vphi_2)$.  For that, it is enough to show that for all $x_i\in S(E,\vphi_i)$ there exists $\xi_i\in S(\mcl{E},\vphi_i)$ such that $\norm{\xi_1-\xi_2}\leq\norm{x_1-x_2}$. Take $\xi_1'\in S(\mcl{E}_0,\vphi_1)$. Let  $U:\sspan\mcl{A}\xi_1'\mcl{B}\rightarrow\sspan\mcl{A} x_1\mcl{B}$ be the bilinear unitary satisfying $U(a\xi_1'b)=ax_1b$. Let $P$ be the bilinear projection of $E$ onto $\sspan\mcl{A} x_1\mcl{B}$. Set
  \begin{align*}
     x_2' &:=Px_2\in\sspan\mcl{A} x_1\mcl{B}\subseteq E, \\
     x_2''&:=(1-P)x_2\in(\sspan\mcl{A} x_1\mcl{B})^\perp\subseteq E, \\
     \vphi_2'(\cdot) &:=\ip{x_2',(\cdot) x_2'}\textnormal{ and } \\
     \vphi_2''(\cdot)&:=\ip{x_2'',(\cdot) x_2''}.
  \end{align*}
  Clearly $\vphi_2=\vphi_2'+\vphi_2''$. Let $\xi_2'=U^*(x_2')\in \sspan\mcl{A} \xi_1'\mcl{B}\subseteq\mcl{E}_0$. Then
  \begin{equation*}
     \ip{\xi_2',a\xi_2'} =\ip{U^*x_2',aU^*x_2'}
                         =\ip{U^*x_2',U^*(ax_2')}
                         =\ip{x_2',ax_2'}
                         =\vphi_2'(a).
  \end{equation*}
  Let $\xi_2''\in S(\mcl{E}_0,\vphi_2'')$. Set $\xi_1=\xi_1'\oplus 0$ and $\xi_2=\xi_2'\oplus\xi_2''$. Then $\xi_i\in S(\mcl{E},\vphi_i)$ with
  \begin{align*}
      \norm{\xi_1-\xi_2}^2 &=\norm{\ip{\xi_1,\xi_1}+\ip{\xi_2,\xi_2}-2\Re{\ip{\xi_1,\xi_2}}} \\
                           &=\norm{\ip{\xi_1',\xi_1'}+\ip{\xi_2',\xi_2'}+\ip{\xi_2'',\xi_2''}-2\Re{\ip{\xi_1',\xi_2'}}} \\
                           &=\norm{\ip{\xi_1'-\xi_2',\xi_1'-\xi_2'}+\ip{\xi_2'',\xi_2''}} \\
                           &=\norm{\ip{U(\xi_1'-\xi_2'),U(\xi_1'-\xi_2')}+\ip{\xi_2'',\xi_2''}} \\
                           &=\norm{\ip{x_1-x_2',x_1-x_2'}+\ip{x_2'',x_2''}} \\
                           &=\norm{\ip{x_1,x_1}+\ip{x_2',x_2'}-2\Re{\ip{x_1,x_2'}}+\ip{x_2'',x_2''}} \\
                           &=\norm{\ip{x_1,x_1}+\ip{x_2,Px_2}-2\Re{\ip{x_1,x_2'}}+\ip{x_2,(1-P)x_2}} \\
                           &=\norm{\ip{x_1,x_1}+\ip{x_2,x_2}-2\Re{\ip{x_1,x_2'}}} \\
                           &=\norm{\ip{x_1-x_2,x_1-x_2}}\qquad\qquad(x_1=x_1\oplus 0, x_2=x_2'\oplus x_2''\textnormal{ in } E) \\
                           &=\norm{x_1-x_2}^2.
  \end{align*}
  Since $x_1,x_2$ and $E$ are arbitrary $\beta_{\mcl{E}}(\vphi_1,\vphi_2)\leq \beta(\vphi_1,\vphi_2)$.

  (ii) Note that $\xi_1\in S(\mcl{E},\vphi_1)$ is independent of $E$ and $\vphi_2$. If we denote $\xi_2$ obtained in part(i) by $\xi_2(x_1,x_2)$, then
  \begin{align*}
      \beta_{\mcl{E}}(\vphi_1,\vphi_2)
                 &=\inf\big\{\norm{\xi-\xi'} : \xi\in S(\mcl{E},\vphi_1),\xi'\in S(\mcl{E},\vphi_2)\big\} \\
                 &\leq\inf\big\{\norm{\xi_1-\xi'} : \xi'\in S(\mcl{E},\vphi_2)\big\} \\
                 &\leq\inf\big\{\norm{\xi_1-\xi_2(x_1,x_2)} : x_i\in S(E,\vphi_i)\big\} \\
                 &=\inf\big\{\norm{x_1-x_2} : x_i\in S(E,\vphi_i)\big\} \\
                 &=\beta_{E}(\vphi_1,\vphi_2)
  \end{align*}
  Since this is true for all common representation module $E$, we get
  \begin{equation*}
       \beta(\vphi_1,\vphi_2) \leq\beta_{\mcl{E}}(\vphi_1,\vphi_2)
                              \leq\inf\big\{\norm{\xi_1-\xi'} : \xi'\in S(\mcl{E},\vphi_2)\big\}
                              \leq\beta(\vphi_1,\vphi_2).
  \end{equation*}
  %
  This completes the proof.
 \eproof

 \bthm\label{betametric}
     $\beta $ is a metric on $CP(\mcl{A},\mcl{B})$.
 \ethm

 \bproof
  \it{Positive definiteness:} Let $\vphi_1,\vphi_2\in CP(\mcl{A},\mcl{B})$. Take $\mcl{E}$ and $\xi_1\in S(\mcl{E},\vphi_1)$ as in proposition \ref{universalbeta}(ii). By definition $\beta(\vphi_1,\vphi_2)\geq 0$. Now if $\beta(\vphi_1,\vphi_2)=0$, then
  \begin{equation*}
       \inf\big\{\norm{\xi_1-\xi_2} : \xi_2\in S(\mcl{E},\vphi_2)\big\}=0.
  \end{equation*}
  Since $S(\mcl{E},\vphi_2)$ is a norm closed subset of $\mcl{E}$, above equality implies that $\xi_1\in S(\mcl{E},\vphi_2)$. Therefore $\vphi_1=\vphi_2$. \\
  \it{Symmetry:} Clear from the definition. \\
  \it{Triangle inequality:} Let $\vphi_1,\vphi_2,\vphi_3\in CP(\mcl{A},\mcl{B})$. Suppose $\mcl{E}$ and $\xi_1\in S(\mcl{E},\vphi_1)$ are as in proposition \ref{universalbeta}(ii). Then
  \begin{align*}
       \beta(\vphi_2,\vphi_3)
              &=\inf\big\{\norm{\xi_2-\xi_3} : \xi_i\in S(\mcl{E},\vphi_i), i=2,3\big\} \\
              &\leq\inf\big\{\norm{\xi_2-\xi_1} : \xi_2\in S(\mcl{E},\vphi_2)\big\}+ \inf\big\{\norm{\xi_1-\xi_3} : \xi_3\in S(\mcl{E},\vphi_3)\big\} \\
              &=\beta(\vphi_2,\vphi_1)+\beta(\vphi_1,\vphi_3).
  \end{align*}
 \eproof

 \subsection{Intertwiners and computation of Bures distance}
  The definition of Bures distance is abstract and does not give us indications as to how to compute it for concrete examples. In this Section, motivated by the work of \cite{KDR08}, we show that Bures distance can be computed using intertwiners between two (minimal) GNS-constructions of CP-maps.

  Suppose $E$ is a common representation module for $\vphi_i$ and $x_i\in S(E,\vphi_i),~i=1,2$. Then $\norm{x_1-x_2} ^2=\norm{\ip{x_1-x_2,x_1-x_2}}=\norm{\vphi_1(1)+\vphi_2(1)-2\Re{\ip{x_1,x_2}}}$. Thus $\beta(\vphi_1,\vphi_2)$ is completely determined by the subsets $\{\ip{x_1,x_2} : x_i\in S(E,\vphi_i)\}\subseteq\mcl{B}$. This observation leads to the following Definition.

 \bdefn
  Given a common representation module $E$  for $\vphi_1$ and $\vphi_2$ define
  \begin{equation*}
       N_E(\vphi_1,\vphi_2):=\big\{\ip{x_1,x_2} : x_i\in S(E,\vphi_i)\big\}
  \end{equation*}
  and
  \begin{equation*}
       N(\vphi_1,\vphi_2):=\underset{E}\cup N_E(\vphi_1,\vphi_2)
  \end{equation*}
  where the union is taken over all common representation module $E$.
 \edefn

  Note that $N(\vphi_1,\vphi_2)\subseteq\mcl{B}$ is always nonempty. Also if $E$ is a common representation module for $\vphi_1$ and $\vphi_2$, then
  \begin{equation}\label{betaE}
      \beta_E(\vphi_1,\vphi_2)=\underset{N\in N_E(\vphi_1,\vphi_2)}{\inf}\,\norm{\vphi_1(1)+\vphi_2(1)-2\Re{N}}^{\frac{1}{2}}
  \end{equation}
  with $\vphi_1(1)+\vphi_2(1)-2\Re{N}=\ip{x_1-x_2,x_1-x_2}\geq 0$.

 \bdefn
  Let $(E_i, x_i)$ be a GNS-construction for $\vphi_i$, $i=1,2$. Then define
  \begin{equation*}
        M(\vphi_1,\vphi_2):=\big\{\ip{x_1,\Phi x_2} : \Phi\in \BL{E_2,E_1},\norm{\Phi}\leq 1\big\}.
  \end{equation*}
 \edefn

 \blem \label{intertwine}
  The set $M(\vphi_1,\vphi_2)\subseteq\mcl{B}$ depends only on the CP-maps $\vphi_i$ and  not on the GNS-constructions $(E_i, x_i)$.
 \elem

 \bproof
  We show that $M(\vphi_1,\vphi_2)$ defined via $(E_i,x_i)$ coincides with $\hat{M}(\vphi_1,\vphi_2)$ which is  defined via the minimal GNS-construction $(\hat{E_i},\hat{x_i})$. Let $U_i:\hat{E_i}\rightarrow \sspan\mcl{A}x_i\mcl{B}$ be the bilinear unitary satisfying $U_i(a\hat{x_i}b)=ax_ib$ for all $a\in\mcl{A},b\in\mcl{B}$. Since $\sspan\mcl{A}x_i\mcl{B}\subseteq E_i$ is a complemented $\mcl{B}$-submodule, $U_i\in \BL{\hat{E_i}, E_i}$ is an adjointable isometry (\cite{lance}, Theorem 3.6). Note that $U_i(\hat{x_i})=x_i$ and $U_i^*(x_i)=\hat{x_i}$. Now suppose $\ip{x_1,\Phi x_2}\in M(\vphi_1,\vphi_2)$, where $\Phi\in\BL{E_2,E_1}$ with $\norm{\Phi}\leq 1$. Set $\hat{\Phi}=U_1^*\Phi U_2$. Then $\hat{\Phi}\in\BL{\hat{E_2},\hat{E_1}}$ with $\snorm{\hat{\Phi}}\leq 1$. Also
  \begin{align*}
       \ip{x_1,\Phi x_2} = \ip{U_1\hat{x_1},\Phi U_2\hat{x_2}}
                         = \ip{\hat{x_1},U_1^*\Phi U_2\hat{x_2}}
                         = \ip{\hat{x_1},\hat{\Phi}\hat{x_2}}
                         \in\hat{M}(\vphi_1,\vphi_2).
  \end{align*}
  Hence $M(\vphi_1,\vphi_2)\subseteq\hat{M}(\vphi_1,\vphi_2)$. To get the reverse inclusion start with a $\hat{\Phi}\in\BL{\hat{E_2},\hat{E_1}}$ and set $\Phi=U_1\hat{\Phi}U_2^*\in\BL{E_2,E_1}$.
 \eproof

 \bprop\label{MN}
  If $(E_i,x_i)$ is a GNS-construction for $\vphi_i,~i=1,2$, then
  \begin{enumerate}
   \item [(i)] $M(\vphi_1,\vphi_2)=N(\vphi_1,\vphi_2)=N_{E_1\oplus E_2}(\vphi_1,\vphi_2)$ and
   \item [(ii)] $\beta(\vphi_1,\vphi_2)
            =\underset{M\in M(\vphi_1,\vphi_2)}{\inf}\,\norm{\vphi_1(1)+\vphi_2(1)-2\Re{M}}^{\frac{1}{2}}.$
   \end{enumerate}
 \eprop

 \bproof
  (i) 
  Suppose $E$ is a common representation module and $\ip{z_1,z_2}\in N_E(\vphi_1,\vphi_2)$. Set $E_1=E_2=E$ and $\Phi=\textnormal{id}_E$. Then, from above Lemma, $\ip{z_1,z_2}=\ip{z_1,\Phi z_2}\in M(\vphi_1,\vphi_2)$. Since $z_1,z_2$ and $E$ are arbitrary $N(\vphi_1,\vphi_2)\subseteq M(\vphi_1,\vphi_2)$. In particular, $M(\vphi_1,\vphi_2)$ is nonempty. For the reverse inclusion, let $\ip{x_1,\Phi x_2}\in M(\vphi_1,\vphi_2)$. Set $z_1=x_1\oplus 0$ and $z_2=\Phi x_2\oplus\sqrt{\textnormal{id}_E-\Phi^*\Phi}x_2$ in $E_1\oplus E_2$. Then $\ip{z_1,az_1}=\ip{x_1,ax_1}=\vphi_1(a)$ and
  \begin{align*}
      \ip{z_2,az_2}
                   &=\ip{\Phi x_2\oplus\sqrt{\textnormal{id}_{E_2}-\Phi^*\Phi}x_2,                                                   a(\Phi x_2)\oplus a\sqrt{\textnormal{id}_{E_2}-\Phi^*\Phi}x_2} \\
                   &=\ip{\Phi x_2,\Phi(ax_2)}+\ip{\sqrt{\textnormal{id}_{E_2}-\Phi^*\Phi}x_2,                                                   \sqrt{\textnormal{id}_{E_2}-\Phi^*\Phi}ax_2} \\
                   &=\ip{x_2,\Phi^* \Phi(ax_2)}+\ip{x_2,(\textnormal{id}_{E_2}-\Phi^*\Phi)ax_2} \\
                   &=\ip{x_2,ax_2} \\
                   &=\vphi_2(a)
  \end{align*}
  for all $a\in\mcl{A}$. Thus $(E_1\oplus E_2,z_i)$ is a GNS-construction for $\vphi_i$. Note that $\ip{x_1,\Phi x_2}=\ip{z_1,z_2}\in N_{E_1\oplus E_2}(\vphi_1,\vphi_2)$. Hence $M(\vphi_1,\vphi_2)\subseteq N_{E_1\oplus E_2}(\vphi_1,\vphi_2)$. Thus $N(\vphi_1,\vphi_2)\subseteq M(\vphi_1,\vphi_2)\subseteq N_{E_1\oplus E_2}(\vphi_1,\vphi_2)\subseteq N(\vphi_1,\vphi_2)$.

  (ii) Follows from equation ($\ref{betaE}$).
 \eproof

 \bcor\label{betainfmod}
  If $(E_i,x_i)$ is GNS-construction for $\vphi_i,~i=1,2$, then
  \begin{equation*}
      \beta(\vphi_1,\vphi_2) =\beta_{E_1\oplus E_2}(\vphi_1,\vphi_2)
                             =\inf\big\{\norm{x_1\oplus 0-y_1\oplus y_2} : y_1\oplus y_2\in S(E_1\oplus E_2,\vphi_2)\big\}.
  \end{equation*}
 \ecor

 \bproof
  Suppose $\ip{x_1,\Phi x_2}\in M(\vphi_1,\vphi_2)$. Then, from the proof of Proposition \ref{MN}, we have $\ip{x_1,\Phi x_2}=\ip{z_1,z_2}$, where $z_i\in S(E_1\oplus E_2,\vphi_i)$ with $z_1= x_1\oplus 0$. Denote the $z_2$ obtained by $z_2(\Phi)$. Then, from proposition \ref{MN}(ii),
  \begin{align*}
      \beta(\vphi_1,\vphi_2)
           &= \inf\big\{\norm{\vphi_1(1)+\vphi_2(1)-2\Re{\ip{x_1\oplus 0,z_2(\Phi)}}}^{\frac{1}{2}} : \Phi\in \BL{E_2,E_1},\norm{\Phi}\leq 1\big\} \\
           &\geq\inf\big\{\norm{\vphi_1(1)+\vphi_2(1)-2\Re{\ip{x_1\oplus 0,y_1\oplus y_2}}}^{\frac{1}{2}} : y_1\oplus y_2\in S(E_1\oplus E_2,\vphi_2)\big\} \\
           &= \inf\big\{\norm{x_1\oplus 0-y_1\oplus y_2} : y_1\oplus y_2\in S(E_1\oplus E_2,\vphi_2)\big\} \\
           &\geq\beta_{E_1\oplus E_2}(\vphi_1,\vphi_2).
  \end{align*}
 \eproof

 \beg
  Let $(X,\mbb{F},\mu)$ be a measure space and let $\mcl{A}=L^{\infty}(X,\mu)$. Consider the states $\vphi_i:\mcl{A}\rightarrow \mbb{C}$ given by
      $\vphi_i(f)=\int fd\mu_i$,
  where $\mu_1$ and $\mu_2$ are two equivalent (i.e., absolutely continuous each other) probability measures on $(X,\mbb{F})$ such that $\mu_i<<\mu, i=1,2$. Let $h$ be a positive function (Radon Nikodym derivative) on $X$ such that $d\mu_1=hd\mu_2$. Clearly $E_i=L^2(X,\mu_i)$ is a von Neumann $\mcl{A}$-$\mbb{C}$-module with left multiplication as the left action. Also  $(E_i,1)$ is a GNS-construction for $\vphi_i$. Suppose $g_1\oplus g_2\in S(E_1\oplus E_2,\vphi_2)$. Then
  \begin{align*}
     \int fd\mu_2 &=\ip{g_1\oplus g_2,f(g_1\oplus g_2)}  \\
                  &=\int |g_1|^2fd\mu_1+\int |g_2|^2fd\mu_2  \\
                  &= \int (|g_1|^2h+|g_2|^2)fd\mu_2
  \end{align*}
  for all $f\in\mcl{A}$, and hence $|g_1|^2h+|g_2|^2=1$ a.e., $\mu_2$.  Therefore
  \begin{align*}
      \beta(\vphi_1,\vphi_2)
                             &=\inf\big\{\norm{1\oplus 0-g_1\oplus g_2} : g_1\oplus g_2\in S(E_1\oplus E_2,\vphi_2)\big\} \\
                             &=\inf\big\{(\ip{1-g_1,1-g_1}+\ip{g_2,g_2})^{\frac{1}{2}} : |g_1|^2h+|g_2|^2=1 \textnormal{ a.e., } \mu_2\big\} \\
                             &=\inf\big\{(2-2\Re{\int g_1d\mu_1})^{\frac{1}{2}} : |g_1|^2h\leq 1 \textnormal{ a.e., } \mu_2\big\} \\
                             &=\sqrt{2}\inf\big\{(1-\int g_1hd\mu_2)^{\frac{1}{2}} : g_1\geq 0\textnormal{ and }0\leq g_1^2h\leq 1 \textnormal{ a.e., } \mu_2\big\} \\
                             &=\sqrt{2}(1-\int \sqrt{h}d\mu_2)^{\frac{1}{2}}.
  \end{align*}
  In particular, if we take $X=\{1,2,\dots,n\}, \mu$ the counting measure, $\mu_1(i)=p_i$ and $\mu_2(i)=q_i$, where $0<p_i,q_i<1$ such that $\sum p_i=\sum q_i=1$, then
  $\beta(\vphi_1,\vphi_2)=\sqrt{2}(1-\sum\sqrt{p_iq_i})^{\frac{1}{2}}$.
 \eeg

  Here we compute the Bures distance for homomorphisms and for some other special cases.

 \bcor\label{betahomo}
  Let $\vphi_1,\vphi_2:\mcl{A}\rightarrow\mcl{B}$ be two unital $\ast$-homomorphisms.
  \begin{enumerate}
    \item [(i)] Then $\beta(\vphi_1,\vphi_2)=\sqrt{2}\inf\big\{\norm{1-\Re{b}}^{\frac{1}{2}} : b\in \mcl{B},\norm{b}\leq 1, \vphi_1(a)b=b\vphi_2(a) \;\forall a\in\mcl{A}\big\}$.
    \item [(ii)] If $\mcl{A}=\mcl{B}$ and $\vphi_2(a)=u^*\vphi_1(a)u$ for some unitary $u\in \mcl{B}$,  then
       \begin{equation*}
         \beta(\vphi_1,\vphi_2)=\sqrt{2}\inf\big\{\norm{1-\Re{b'u}}^{\frac{1}{2}}:b'\in\vphi_1(\mcl{A})', \norm{b'}\leq 1\big\}.
       \end{equation*}
    \item [(iii)] If $u\in M_n(\mbb{C})$ is a unitary and  $\vphi: M_n(\mbb{C})\rightarrow M_n(\mbb{C})$ is the $\ast$-homomorphism $\vphi(a)=u^*au$, then
       \begin{equation*}
        \beta(id.,\vphi)=\sqrt{2}\inf\big\{\norm{1-\Re{\lambda u}}^{\frac{1}{2}} : \lambda\in [-1,1]\big\}.
       \end{equation*}
  \end{enumerate}
 \ecor

 \bproof
  (i) Let $E_i$ be the von Neumann $\mcl{A}$-$\mcl{B}$-module $\mcl{B}$ with left action $ax:=\vphi_i(a)x$ for all $a\in\mcl{A},~x\in E_i$. Then $(E_i,1)$ is the minimal GNS-construction for $\vphi_i$.  Suppose $\Phi\in\BL{E_2,E_1}$. Then 
  \begin{equation*}
      \vphi_1(a)\Phi(1)=a\Phi(1)=\Phi(a1)=\Phi(\vphi_2(a))=\Phi(1)\vphi_2(a)
  \end{equation*}
  for all $a\in\mcl{A}$. Clearly, for a fixed $b_0\in\mcl{B}$ satisfying $\vphi_1(a)b_0=b_0\vphi_2(a)$, the map $b\mapsto b_0b$ is an element of  $\BL{E_2,E_1}$. Thus
  \begin{align*}
      \beta(\vphi_1,\vphi_2)
              &=\inf\big\{\norm{\vphi_1(1)+\vphi_2(1)-2\Re{M}}^{\frac{1}{2}} : M\in M(\vphi_1,\vphi_2)\big\} \\
              &=\inf\big\{\norm{2-2\Re{\Phi(1)}}^{\frac{1}{2}}:\Phi\in\BL{E_2,E_1} ,\norm{\Phi}\leq 1 \big\} \\
              &=\sqrt{2}\inf\big\{\norm{1-\Re{b}}^{\frac{1}{2}} : b\in \mcl{B}, \norm{b}\leq 1,  \vphi_1(a)b=b\vphi_2(a)\;\forall a\in\mcl{A}\big\}.
  \end{align*}

  (ii) Suppose $b\in \mcl{B}$. Then $\vphi_1(a)b=b\vphi_2(a)$ for all $a\in\mcl{A}$ implies that $bu^* \in \vphi_1(\mcl{A})'$, and hence $b=b'u$ for some $b'\in\vphi_1(\mcl{A})'\subseteq \mcl{B}$.\\

  (iii) This follows from (ii), since $M_n'=\mbb{C}I$.
 \eproof

  In \cite{KDR08} it is shown that the Bures distance is comparable with completely bounded norm when $\mcl{B}=\B{G}$, and the following bounds were obtained. In fact, the lower bound holds even for an  arbitrary unital $C^*$-algebra $\mcl{B}$.

%

 \bthm[\cite{KDR08}]\label{KDRthm}
  For $\vphi_1,\vphi_2\in CP(\mcl{A},\B{G}),$
  \begin{equation*}
       \dfrac{\norm{\vphi_1-\vphi_2}_{cb}}{\sqrt{\norm{\vphi_1}_{cb}}+\sqrt{\norm{\vphi_2}_{cb}}}\leq \beta(\vphi_1,\vphi_2)\leq\sqrt{\norm{\vphi_1-\vphi_2}_{cb}}.
  \end{equation*}
  Moreover, there exists a common representation module $E$ and corresponding GNS-construction $(E,x_i)$ for $\vphi_i$ such that
  \begin{equation*}
      \beta(\vphi_1,\vphi_2)=\beta_E(\vphi_1,\vphi_2)=\norm{x_1-x_2}.
  \end{equation*}
 \ethm

 \beg\label{betaopnormeg}
  In general, the upper bound given in Theorem \ref{KDRthm} may fails to hold if the cb-norm is replaced by the operator norm. For example, consider the CP-maps $\vphi_i: M_2(\mbb{C})\rightarrow M_2(\mbb{C})$ given by
  \begin{equation*}
      \vphi_1(\begin{bmatrix}
                       a_{ij}
             \end{bmatrix}):=
             \begin{bmatrix}
                   a_{11}+2a_{22} & a_{21} \\
                   a_{12} & a_{22}+2a_{11}
             \end{bmatrix}
            \quad\textnormal{ and }\quad
     \vphi_2(\begin{bmatrix}
                     a_{ij}
             \end{bmatrix}):=
             \begin{bmatrix}
                    2a_{22} & 0 \\
                    0 & 2a_{11}
             \end{bmatrix}.
  \end{equation*}
  Let $E=M_{8\times 2}(\mbb{C})$ which is a von Neumann $M_2(\mbb{C})$-$M_2(\mbb{C})$-module with module actions given by
  \begin{equation*}
       axb:=\begin{bmatrix}
                   ax_1b \\
                   ax_2b  \\
                   ax_3b \\
                   ax_4b
            \end{bmatrix}\qquad\forall x=\begin{bmatrix}
                                                 x_1\\
                                                 x_2\\
                                                 x_3\\
                                                 x_4
                                         \end{bmatrix}\in E
                                         \textnormal{ and } a,b,x_i\in M_2(\mathbb{C}).
  \end{equation*}
  Then $E$ is a common representation module with
  \begin{align*}
   z_1:&=\left[\begin{array}{cccccccc}
              1 & 0 & 0 & 0 & 0 & \frac{\sqrt{3}}{\sqrt{2}} & 0 & \frac{-1}{\sqrt{2}} \\
              0 & 0 & 0 & 1 & \frac{\sqrt{3}}{\sqrt{2}} & 0 & \frac{1}{\sqrt{2}} & 0
              \end{array}
        \right]^{\textnormal{t}} \in S(E,\vphi_1)
   \end{align*}
   and
   \begin{align*}
   z_2:&=\left[\begin{array}{cccccccc}
              0 & 0 & 0 & 0 & 0 & 1 & 0 & -1 \\
              0 & 0 & 0 & 0 & 1 & 0 & 1 & 0
              \end{array}
        \right]^{\textnormal{t}}\in S(E,\vphi_2),
  \end{align*}
  where `$\textnormal{t}$' stands for transpose.
%
  Note that if $x\oplus y=[x_{ij}]\oplus[y_{ij}]\in S(E\oplus E,\vphi_2)$, then by evaluating $\vphi _2$ at matrix units, we see that $x_{i1}=y_{i1}=0=x_{k2}=y_{k2},~i=1,3,5,7,~ k=2,4,6,8$ and
  \begin{equation*}
   \left.\begin{aligned}
           &\sum_{i=2,4,6,8}(\ol{x_{i1}}x_{i-1,2}+\ol{y_{i1}}y_{i-1,2})=0, \\
           &\sum_{i=2,4,6,8}(\abs{x_{i1}}^2+\abs{y_{i1}}^2)=2=\sum_{i=1,3,5,7}(\abs{x_{i2}}^2+\abs{y_{i2}}^2).
        \end{aligned}
        \qquad
   \right\} \qquad (*)
  \end{equation*}
  Hence
  \begin{align*}
     \beta(\vphi_1,\vphi_2)
         &=\inf\,\norm{z_1\oplus 0-x\oplus y} \\
         &=\inf\,\norm{\begin{bmatrix}
                           \medskip
                           5-\Re{\sqrt{6}x_{61}-\sqrt{2}x_{81}} & -x_{12}-\ol{x_{41}} \\
                           -\ol{x_{12}}-x_{41}   & 5-\Re{\sqrt{6}x_{52}+\sqrt{2}x_{72}}
                        \end{bmatrix}
                      }^{\frac{1}{2}} \\
         &\geq \inf\,\norm{\begin{bmatrix}
                                \medskip
                                5-\Re{\sqrt{6}x_{61}-\sqrt{2}x_{81}} & 0 \\
                                0 &  5-\Re{\sqrt{6}x_{52}+\sqrt{2}x_{72}}
                            \end{bmatrix}
                          }^{\frac{1}{2}}
  \end{align*}
  where the infimums are taken over all $x\oplus y\in E_1\oplus E_2$ satisfying $(*)$.  Now some elementary calculus shows that $\beta (\vphi_1,\vphi_2)\geq \sqrt{5-\sqrt{2}-\sqrt{6}}.$
  Note that $\norm{z_1-z_2}=\sqrt{5-\sqrt{2}-\sqrt{6}}$, and hence $\beta(\vphi_1,\vphi_2)=\sqrt{5-\sqrt{2}-\sqrt{6}}>1$. But $\vphi_1-\vphi_2$ is the transpose map. Therefore $1=\norm{\vphi_1-\vphi_2}<\beta(\vphi_1,\vphi_2)^2 <\norm{\vphi_1-\vphi_2}_{cb}=2$ (see \cite {paulsen} for the computation of cb-norm for transpose map).
 \eeg

  Theorem \ref{KDRthm} guarantees the existence of a common representation module, where Bures distance is attained. It is a natural question as to whether Bures distance is attained in every common representation module. This is true for states (\cite{araki}). The question in the general case was asked by \cite{KDR06,KDR08}. Here we resolve it in the negative through a simple counter example.

 \beg\label{betamineg}
  Consider the (normal) CP-maps $\vphi_i: M_2(\mbb{C})\rightarrow M_2(\mbb{C})$ given by $\vphi_i(a):=a_i^*aa_i$, where
     $a_1=\begin{bmatrix}
                1 & 0 \\
                0 & 0
          \end{bmatrix}$ and
     $a_2=\begin{bmatrix}
                0 & 1 \\
                0 & 0
         \end{bmatrix}$.
  Then $(\hat{E_i},\hat{x_i}):=(M_2(\mbb{C}),a_i)$ is the minimal GNS-construction for $\vphi_i$. Set $x_1=\hat{x}_1\oplus 0$ and $x_2=0\oplus\hat{x}_2$. Then $x_i\in S(\hat{E_1}\oplus \hat{E_2},\vphi_i)$ and
  \begin{align*}
       \beta(\vphi_1,\vphi_2)=\beta_{\hat{E_1}\oplus \hat{E_2}}(\vphi_1,\vphi_2)
                             \leq \norm{x_1-x_2}
                             =\norm{I}
                             =1.
  \end{align*}
  Clearly, $E:=M_2(\mbb{C})$ is a common representation module. It is not hard to see that $S(E,\vphi_i) =\{\lambda a_i : \lambda\in \mbb{C},\abs{\lambda}= 1\}$. Now for any $x_i=\lambda_i a_i\in S(E,\vphi_i)$,
  \begin{align*}
      \norm{x_1-x_2}^2 =\norm{\begin{bmatrix}
                                     \medskip
                                     1  &  -\ol{\lambda}_1\lambda_2 \\
                                    -\ol{\lambda}_2\lambda_1  & 1
                              \end{bmatrix}}
                       =\sup\,\Bigg\{\abs{\lambda} : \lambda\in\sigma\Big(\begin{bmatrix}
                                                                        \medskip
                                                                        1  &  -\ol{\lambda}_1\lambda_2 \\
                                                                        -\ol{\lambda}_2\lambda_1  & 1                                                              \end{bmatrix}\Big)\Bigg\}
                       =2.
  \end{align*}
  Hence $\beta_E(\vphi_1,\vphi_2)=\sqrt{2}>1\geq \beta(\vphi_1,\vphi_2)$. Note that here $\beta(\vphi_1,\vphi_2)\leq 1=\sqrt{\norm{\vphi_1-\vphi_2}}$.
 \eeg

 \bconj
  If $\vphi,\psi\in CP(\mcl{A},\mcl{B})$, then
  \begin{equation*}
      \beta(\vphi,\psi)=\underset{\phi,n}{\sup}\,\beta(\phi\circ\vphi_n,\phi\circ\psi_n)
  \end{equation*}
  where the supremum is taken over all states $\phi:M_n(\mcl{B})\rightarrow\mbb{C}, n\geq 1$.
 \econj

  From Proposition \ref{betaampl} and \ref{betacomp}  we have $\beta(\phi\circ\vphi_n,\phi\circ\psi_n)\leq \beta(\vphi_n,\psi_n)=\beta(\vphi,\psi)$ for all states $\phi:M_n(\mcl{B})\rightarrow\mbb{C}, n\geq 1$. If the
  conjecture can be proved directly, then using the upper bound for states \cite{Bures,KDR08} we get an
  alternative proof of the upper bound for Bures metric:
  \begin{equation*}
      \beta(\vphi,\psi)=\underset{\phi,n}{\sup}\,\beta(\phi\circ\vphi_n,\phi\circ\psi_n)
                       \leq \underset{\phi,n}{\sup}\,\sqrt{\norm{\phi\circ\vphi_n-\phi\circ\psi_n}}
                       = \sqrt{\norm{\varphi-\psi}_{cb}}.
  \end{equation*}

 \section{Bures distance: \texorpdfstring{$C^*$}{C*}-algebras}
  This Section consists mostly of counter examples. But results similar to the last section do hold for injective $C^*$-algebras.

 \subsection{Counter examples}
  We saw that if the range algebras are von Neumann algebras, then the Bures metric can be computed using intertwiners. It was crucial that the space of intertwiners was independent of the choice of GNS-constructions (Lemma \ref{intertwine} ). The first example here shows that this is no longer the case for some range $C^*$-algebras. We have another example to show that the upper bound computed for $\beta $ in Theorem \ref{KDRthm} may not hold for general range $C^*$-algebras. Finally, as a worst case scenario we have a tricky example to show that even the triangle inequality may fail to hold.

 \beg\label{intertwine-eg}
  If $\vphi_1$ and $\vphi_2$ are CP-maps between $C^*$-algebras, then $M(\vphi_1,\vphi_2)$ may depends on the GNS-construction. For example, consider the CP-maps $\vphi_i:C([0,2\pi])\rightarrow C([0,2\pi])$  given by $\vphi_i(f):=g_if$, where $g_i(t)=|sin(t)|^i$ for all $t\in [0,2\pi],~i=1,2$. Set $\hat{x_i}=\sqrt{g_i}$ and
  \begin{align*}
      \hat{E_i}&=\cspan\,\{\sqrt{g_i}f: f\in C([0,2\pi])\} \\
               &=\{f\in C([0,2\pi]): f(0)=f(\pi)=f(2\pi)=0\}.
  \end{align*}
  Then $(\hat{E_i},\hat{x_i})$ is the minimal GNS-construction for $\vphi_i$. Define the adjointable bilinear map $\hat{\Phi}:\hat{E_2}\rightarrow \hat{E_1}$ by $\hat{\Phi}(f)=gf$, where
  \begin{align*}
      g(t)=\left\{\begin{array}{lll}
                       \medskip
                       \frac{1}{2} & \textnormal{if} & 0\leq t < \pi, \\
                       1 & \textnormal{if} & \pi \leq t \leq 2\pi.
                 \end{array}\right.
  \end{align*}
  Since $\hat{\Phi}$ is a contraction $\ip{\hat{x_1},\hat{\Phi}\hat{x_2}}\in \hat{M}(\vphi_1,\vphi_2)$. We have $(E_i,x_i):=(C([0,2\pi]),\hat{x_i})$ is also a GNS-construction for $\vphi_i$. Now if $\Phi:E_2\rightarrow E_1$ is an adjointable bilinear map, then $\Phi(f)=\Phi(1)f$ for all $f\in C([0,2\pi])$. Thus $\BL{E_2,E_1}=\{f\mapsto hf: h\in C([0,2\pi])\}$. Hence if $\ip{\hat{x_1},g\hat{x_2}}=\ip{\hat{x_1},\hat{\Phi}\hat{x_2}}\in M(\vphi_1,\vphi_2)$, then $\ip{\hat{x_1},g\hat{x_2}}=\ip{\hat{x_1},h\hat{x_2}}$ for some $h\in C([0,2\pi])$; i.e.,
  \begin{equation*}
     \begin{array}{rll}
         \hat{x_1}(t)g(t)\hat{x_2}(t) &=\hat{x_1}(t)h(t)\hat{x_2}(t), &\forall  t\in[0,2\pi] \\
         \Rightarrow \qquad      g(t) &=h(t), &\forall t\in [0,2\pi]\smallsetminus\{0,\pi,2\pi\}
     \end{array}
  \end{equation*}
  which is not possible since $h$ is continuous on $[0,2\pi]$. So $\ip{\hat{x_1},\hat{\Phi}\hat{x_2}}\notin M(\vphi_1,\vphi_2)$.
 \eeg

 \beg\label{betacbnormeg}
  Suppose $H$ is an infinite dimensional Hilbert space and $p\in\B{H}$ is an orthogonal projection such that both $p$ and $q:=(1-p)$ have infinite rank. Let $\mcl{A}=C^*\{\K{H}\cup \{I\}\}$ and let $u=\lambda p+\ol{\lambda}q$, where  $\lambda = e^{i\theta }$ is a scalar  with $-\frac{\pi}{2} <\theta <\frac{\pi}{2}$. Note that $u\in\B{H}$ is a unitary. Define $\ast$-homomorphisms $\vphi_i:\mcl{A}\rightarrow\mcl{A}$ by $\vphi_1(a):=a$ and $\vphi_2(a):=u^*au$. Now suppose $E$ is a common representation module for $\vphi_1,\vphi_2$ and $x_i\in S(E,\vphi_i)$. Since $\norm{ax_i-x_i\vphi_i(a)}=0$, we get $ax_i=x_i\vphi_i(a)$ for all $a\in\mcl{A}$. Then
  \begin{align*}
        a\ip{x_1,x_2}=\vphi_1(a)\ip{x_1,x_2}
                     =\ip{x_1,x_2}\vphi_2(a)
                     =\ip{x_1,x_2}u^*au
  \end{align*}
  for all $a\in\mcl{A}$, and hence $\ip{x_1,x_2}u^*\in \mcl{A}'$. Therefore $\ip{x_1,x_2}=\lambda'u$ for some $\lambda' \in\mbb{C}$. Since $\ip{x_1,x_2}\in \mcl{A}$ and $u\not\in\mcl{A}$ we have $\lambda'=0$, whence $\ip{x_1,x_2}=0$. Also since $E$ and $x_i\in S(E,\vphi_i)$ are arbitrary
  \begin{align*}
       \beta(\vphi_1,\vphi_2)=\underset{E,x_i}{\inf}\,\norm{x_1-x_2}
                             =\norm{\vphi_1(1)+\vphi_2(1)}^{\frac{1}{2}}
                             =\sqrt{2}.
  \end{align*}
  Now we prove that $\sqrt{\norm{\vphi_1-\vphi_2}_{cb}}<\beta(\vphi_1,\vphi_2)$. For $a=[a_{ij}]\in \B{H}=\B{H_p\oplus H_p^\perp}$, where $H_p=\ran{p}$,
  \begin{align*}
     \norm{\vphi_1(a)-\vphi_2(a)}  & = \norm{a-u^*au} \\
                                   &= \norm{ \begin{bmatrix}
                                                 a_{11} & a_{12} \\
                                                 a_{21} & a_{22}
                                             \end{bmatrix}-
                                             \begin{bmatrix}
                                                 \lambda & 0 \\
                                                 0 & \ol{\lambda}
                                             \end{bmatrix}^*
                                             \begin{bmatrix}
                                                 a_{11} & a_{12} \\
                                                 a_{21} & a_{22}
                                             \end{bmatrix}
                                             \begin{bmatrix}
                                                 \lambda & 0 \\
                                                 0 & \ol{\lambda}
                                             \end{bmatrix}
                                           }\\
                                   &= \norm{ \begin{bmatrix}
                                                 0 & (1-\ol{\lambda}^2)a_{12}\\
                                                 (1-\lambda^2)a_{21} & 0
                                             \end{bmatrix}
                                           } \\
                                   &= \max\big\{\bnorm{(1-\ol{\lambda}^2)a_{12}}, \norm{(1-\lambda^2)a_{21}}\big\} \\
                                   & \le\abs{1-\lambda^2}\norm{a}
  \end{align*}
  so that $\norm{\vphi_1-\vphi_2}\leq \abs{1-\lambda^2}$. But $a=\begin{bmatrix}0&I\\I&0\end{bmatrix}$ is of norm one and $\norm{(\vphi_1-\vphi_2)(a)}=\abs{1-\lambda^2}$, whence $\norm{\vphi_1-\vphi_2}=\abs{1-\lambda^2}=\abs{\lambda(\ol{\lambda}-\lambda)}=\abs{\ol{\lambda}-\lambda}$. Now for all $n\geq 1$, if we let $U_n ,P_n$ and $Q_n$ denote the $n\times n$ diagonal matrix with diagonal $u,p$ and $q$ respectively, then  $U_n=\lambda P_n+\ol{\lambda} Q_n$ and $(\vphi_1-\vphi_2)_n(A)=A-U_n^*AU_n$ for all $A\in M_n(\mcl{A})$. Then, as above, we get $\norm{(\vphi_1-\vphi_2)_n}=\abs{\ol{\lambda}-\lambda}$. Thus
  \begin{equation*}
       \sqrt{\norm{\vphi_1-\vphi_2}}=\sqrt{\norm{\vphi_1-\vphi_2}_{cb}}
                                    =\sqrt{\abs{\ol{\lambda}-\lambda}}
                                    <\sqrt{2}
                                    =\beta(\vphi_1,\vphi_2).
  \end{equation*}

  Now if $\vphi_i$ is considered as a map into $\B{H}$ denote it by $\wtilde{\vphi}_i$. Then  $b\in\wtilde{\vphi}_1(\mcl{A})'\subseteq\B{H}$ implies that $ba=ab$ for all $a\in\K{H}\subseteq\mcl{A}$, so that $b=\lambda_b I$ for some $\lambda_b\in\mbb{C}$. From Corollary \ref{betahomo},
  \begin{align*}
       \beta(\wtilde{\vphi}_1,\wtilde{\vphi}_2)
             &=\sqrt{2}\inf\big\{\norm{1-\Re{\lambda' u}}^{\frac{1}{2}} : \lambda'\in\mbb{C}, \abs{\lambda'}\leq 1\big\} \\
             &\leq\sqrt{2}\,\norm{1-\Re{u}}^{\frac{1}{2}} \\
             &=\sqrt{2}\abs{1-\Re{\lambda}}^{\frac{1}{2}} \\
             &< \sqrt{2}\\
             &=\beta(\vphi_1,\vphi_2).
  \end{align*}
 \eeg

 \beg\label{betametriceg}
  Let $H$ be an infinite dimensional Hilbert space. Consider the unital  $C^*$-subalgebra
  \begin{align*}
       \mcl{A}:&=C^*\Bigg\{\K{H\oplus H} \cup \Big\{
                                                  \begin{bmatrix}
                                                         I&0\\
                                                         0&0
                                                  \end{bmatrix},
                                                  \begin{bmatrix}
                                                         0&0\\
                                                         0&I
                                                  \end{bmatrix}
                                             \Big\}
                   \Bigg\} \\
               &=  \Bigg\{
                        \begin{bmatrix}
                              \lambda_1I+a_{11} &   a_{12}\\
                              a_{21}            &  \lambda_2I+a_{22}
                        \end{bmatrix}
                        : \lambda_i\in\mbb{C},\,a_{ij}\in\K{H}
                   \Bigg\}
  \end{align*}
  of $\B{H\oplus H}$, where $\K{\cdot}$ is the set of compact operators. Suppose $u\in\B{H}$ is a unitary and $1<r\in\mbb{R}$. Set
  \begin{equation*}
     z_1=\begin{bmatrix}
                0&u\\
                0&rI
         \end{bmatrix},\,
     z_2=\begin{bmatrix}
                0&0\\
                0&rI
         \end{bmatrix}
     \textnormal{ and }
     z_3=\begin{bmatrix}
                0&I\\
                0&rI
         \end{bmatrix}
  \end{equation*}
  in $\B{H\oplus H}$. Define CP-maps $\vphi_i:\mcl{A}\rightarrow\mcl{A}$ by $\vphi_i(a):=z_i^*az_i,\,i=1,2,3.$ Note that each $\vphi _i$ has the form, $\vphi_i(\cdot)=\begin{bmatrix}0&0\\0&\ast\end{bmatrix}$. Let
  \begin{equation*}
      E_{12}=\Bigg\{
                   \begin{bmatrix}
                          x_{11} &  \lambda_1u+x_{12}\\
                          x_{21} &  \lambda_2I+x_{22}
                   \end{bmatrix}
                   : \lambda_i\in \mbb{C},\,x_{ij}\in\K{H}
              \Bigg\}
  \end{equation*}
  which is a Hilbert $\mcl{A}$-$\mcl{A}$-module with a natural inner product and bimodule structure. Note that $z_i\in S(E_{12},\vphi_i),\,i=1,2$, and hence $\beta(\vphi_1,\vphi_2)\leq\norm{z_1-z_2}=1$. Similarly
  \begin{equation*}
      E_{23}=\Bigg\{
                   \begin{bmatrix}
                         x_{11} &  \lambda_1I+x_{12}\\
                         x_{21} &  \lambda_2I+x_{22}
                   \end{bmatrix}
                   : \lambda_i\in \mbb{C},\,x_{ij}\in\K{H}
              \Bigg\}
  \end{equation*}
  is a Hilbert $\mcl{A}$-$\mcl{A}$-module with $z_i\in S(E_{23},\vphi_i),\,i=2,3$, and $\beta(\vphi_2,\vphi_3)\leq\norm{z_2-z_3}=1$. Now we will show that $\beta(\vphi_1,\vphi_3)> 2\geq\beta(\vphi_1,\vphi_2)+\beta(\vphi_2,\vphi_3)$ so that $\beta$ fails to satisfy triangle inequality. Suppose $E$ is a common representation module for $\vphi_1,\vphi_3$. We prove that $\ip{x_1,x_3}=0$ for all $x_i\in S(E, \vphi_i)$. If we proved this, then $E$ and $x_i\in S(E,\vphi_i)$  arbitrary implies that
  \begin{equation*}
       \beta(\vphi_1,\vphi_3)=\underset{E,x_i}{\inf}\,\norm{x_1-x_3}
                             =\norm{\vphi_1(1)+\vphi_3(1)}^{\frac{1}{2}}
                             =\sqrt{2(1+r^2)}>2.
  \end{equation*}
  Suppose $\ip{x_1,x_3}=[a_{ij}]$. Since
     $ 0 \leq\begin{bmatrix}
                   \ip{x_1,x_1} & \ip{x_1,x_3} \\
                   \ip{x_3,x_1} & \ip{x_3,x_3}
             \end{bmatrix}
         =\left[\begin{array}{cc|cc}
                             0        &    0      & a_{11} &  a_{12}\\
                             0        &   \ast    & a_{21} &  a_{22} \\ \hline
                             a_{11}^* & a_{21}^*  & 0      &  0 \\
                             a_{12}^* & a_{22}^*  & 0      & \ast
                \end{array}\right]$
  we have $a_{11}=a_{12}=a_{21}=0$. Also for all $a\in\K{H}$, we get
  \begin{equation*}
       \begin{bmatrix}
              a&0\\
              0&0
       \end{bmatrix}x_1
       =
       x_1\begin{bmatrix}
                 0&0\\
                 0&u^*au
          \end{bmatrix}
          \quad\textnormal{and}\quad
          \begin{bmatrix}
                 a&0\\
                 0&0
          \end{bmatrix}x_3
        =
        x_3\begin{bmatrix}
                  0&0\\
                  0&a
           \end{bmatrix}.
  \end{equation*}
  (Simply look at the norm of the difference.) Hence
  \begin{equation*}
       \begin{bmatrix}
               0&0\\
               0&u^*au
       \end{bmatrix}
       \ip{x_1,x_3}
       =
       \ip{x_1,x_3}
       \begin{bmatrix}
               0&0\\
               0&a
       \end{bmatrix};
  \end{equation*}
  i.e.,
  \begin{equation*}
       \begin{bmatrix}
              0&0\\
              0&u^*au
       \end{bmatrix}
       \begin{bmatrix}
              0&0\\
              0&a_{22}
       \end{bmatrix}
       =
       \begin{bmatrix}
              0&0\\
              0&a_{22}
       \end{bmatrix}
       \begin{bmatrix}
              0&0\\
              0&a
       \end{bmatrix}
  \end{equation*}
  which implies that $u^*aua_{22}=a_{22}a$; i.e., $aua_{22}=ua_{22}a$ for all $a\in\K{H}$. Hence  $ua_{22}=\lambda I$ for some $\lambda\in\mbb{C}$. Thus $a_{22}=\lambda u^*$. Since $a_{22}\in\K{H}$ and $u^*\notin\K{H}$ we have $\lambda=0$, and hence $a_{22}=0$ and $\ip{x_1,x_3}=0$.
 \eeg

 \subsection{Injective \texorpdfstring{$C^*$}{C*}-algebras}
  Recall that a $C^*$-algebra $\mcl{B}$ is an injective $C^*$-algebra if, whenever $\mcl{C}$ is a $C^*$-algebra, $\mcl{S}$ an operator system contained in $\mcl{C}$, and $\vphi:\mcl{S}\rightarrow\mcl{B}$ is a completely positive contraction, then $\vphi$ extends to a completely positive contraction $\tilde{\vphi}:\mcl{C}\rightarrow\mcl{B}$. Further, this is equivalent to saying that there is a faithful representation $\pi $ of $\mcl{B}$ on a Hilbert space $G$, such that
  there is a conditional expectation from $\B{G}$ onto $\pi(\mcl{B})$. See \cite{Arv,paulsen,Takesaki} for details.

 \bprop\label{injectivebeta}
  Let $\mcl{A}$ and $\mcl{B}$ be unital $C^*$-algebras and  $\vphi_1,\vphi_2 \in CP(\mcl{A},\mcl{B})$.
  Suppose $\mcl{B}$ is an injective unital $C^*$-algebra and $\pi:\mcl{B}\rightarrow \B{G}$ is a faithful representation of $\mcl{B}$ on $G$. Then $\beta(\vphi_1,\vphi_2)=\beta(\pi\circ\vphi_1,\pi\circ\vphi_2)$.
 \eprop

 \bproof
  Since $\mcl{B}$ is injective there exists a completely  positive conditional expectation $P:\B{G}\rightarrow \pi(\mcl{B})$. Take $\vphi=\pi^{-1}\circ P:\B{G}\rightarrow\mcl{A}$. Then $\vphi$ is a contractive CP-map. Moreover, $\vphi\circ\pi\circ\vphi_i=\vphi_i,~i=1,2.$  Now by Proposition \ref{betacomp},
  \begin{equation*}
    \beta(\vphi_1,\vphi_2)=\beta(\vphi\circ\pi\circ\vphi_1,\vphi\circ\pi\circ\vphi_2)
                          \leq\beta(\pi\circ\vphi_1,\pi\circ\vphi_2)
                          \leq\beta(\vphi_1,\vphi_2).
  \end{equation*}
 \eproof

  From Proposition \ref{betacomp}, we know that $\beta(\pi\circ\vphi_1,\pi\circ\vphi_2)\leq\beta(\vphi_1,\vphi_2)$ even for an arbitrary $C^*$-algebra $\mcl{B}$. But, in general, equality may not holds. See example \ref{betacbnormeg}.

  The following bounds were first obtained in \cite{KDR08}.

 \bcor
  If $\mcl{B}$ is an injective unital $C^*$-algebra, then $\beta$ is a metric on $CP(\mcl{A},\mcl{B})$ and
  \begin{equation*}
     \dfrac{\norm{\vphi_1-\vphi_2}_{cb}}{\sqrt{\norm{\vphi_1}_{cb}}+\sqrt{\norm{\vphi_2}_{cb}}}
            \leq\beta(\vphi_1,\vphi_2)
            \leq\sqrt{\norm{\vphi_1-\vphi_2}_{cb}}.
  \end{equation*}
  Further, there exists common representation module $E$ and corresponding GNS-construction $(E,x_i)$ for $\vphi_i$ such that
  \begin{equation*}
   \beta(\vphi_1,\vphi_2)=\beta_E(\vphi_1,\vphi_2)=\norm{x_1-x_2}.
  \end{equation*}
 \ecor

 \bproof
  Suppose $\pi:\mcl{B}\rightarrow\B{G}$ is a faithful representation of $\mcl{B}$. Now the first part follows from Theorem \ref{betametric} and Proposition \ref{injectivebeta}. 
  Also from Theorem \ref{KDRthm} and Proposition \ref{injectivebeta}, we have
  \begin{equation*}
    \beta(\vphi_1,\vphi_2)
            =\beta(\pi\circ\vphi_1,\pi\circ\vphi_2)
            \leq\sqrt{\norm{\pi\circ\vphi_1-\pi\circ\vphi_2}_{cb}}
            =\sqrt{\norm{\vphi_1-\vphi_2}_{cb}}.
  \end{equation*}
  Now, from Theorem \ref{KDRthm}, we know that there exists a von Neumann $\mcl{A}$-$\B{G}$-module $F$ with $y_i\in S(F,\pi\circ\vphi_i)$ such that $\norm{y_1-y_2}=\beta(\pi\circ\vphi_1,\pi\circ\vphi_2)$.  Given $b\in\mcl{B},~y\in F$ define $yb:=y\pi(b)$. Under this action, $F$ forms a right $\mcl{B}$-module, denoted by $E_0$. Let $P:\B{G}\to\pi(\mcl{B})$ be a completely positive conditional expectation satisfying $P(b_1ab_2)=b_1P(a)b_2$ for all $b_i\in\pi(\mcl{B}),~a\in\B{G}$. Now define a $\mcl{B}$-valued semi-inner product on $E_0$ by $\ip{x_1,x_2}':=\pi^{-1}P(\ip{x_1,x_2})$. Let $E$ be the completion of the $\mcl{B}$-valued inner product space $E_0/N$, where $N:=\{ x\in E_0: \ip{x,x}'=0\}$. Then $E$ is a Hilbert $\mcl{A}$-$\mcl{B}$-module with left action induced by that of $\mcl{A}$ on $F$. Note that $x_i:=y_i+N\in S(E,\vphi_i),~i=1,2$ are such that
  \begin{align*}
      \beta_E(\vphi_1,\vphi_2) \leq\norm{x_1-x_2}
                               \leq\norm{y_1-y_2}
                               =\beta(\pi\circ\vphi_1,\pi\circ\vphi_2)
                               =\beta(\vphi_1,\vphi_2).
  \end{align*}
  Thus $\beta(\vphi_1,\vphi_2)=\beta_E(\vphi_1,\vphi_2)=\norm{x_1-x_2}$.
 \eproof

 \section{Bures distance and a rigidity theorem}
  Observe that for the identity map on a unital $C^*$-algebra $\mcl{B}$ the GNS-module is $\mcl{B}$ itself. Here we show that if a CP-map on a von Neumann algebra $\mcl {B}$ is close to the identity map in Bures distance then the GNS-module has a copy of $\mcl{B}$.

  Suppose $\mcl{B}\subset\B{G}$ is a von Neumann algebra and $\vphi:\mcl{B}\to\mcl{B}$ is a CP-map.

 \bprop\label{domination}
  If $(E,x)$ is the minimal GNS-construction for $\vphi$, then the following are equivalent:
  \begin{enumerate}
     \item [(i)] The center $C_{\mcl{B}}(E):=\{y\in E: by=yb\quad\forall b\in\mcl{B}\}$ contains a unit vector.
     \item [(ii)] $E\cong\mcl{B}\oplus F$ for some von Neumann $\mcl{B}$-$\mcl{B}$-module $F$.
     \item [(iii)] There exists an  element $c\in\mcl{B}$ such that the two sided (strongly closed) ideal generated by $c$ is $\mcl{B}$, and a CP-map $\psi:\mcl{B}\to\mcl{B}$ such that
        \begin{equation*}
           \vphi(b)=c^*bc+\psi(b)\quad\forall b\in\mcl{B}.
        \end{equation*}
  \end{enumerate}
 \eprop

 \bproof
  $(i)\Rightarrow (ii)$: Let $z\in C_{\mcl{B}}(E)$ be a unit vector. The two sided $\mcl{B}$-$\mcl{B}$-module generated by $z$ is naturally isomorphic to $\mcl{B}$ by $bz\mapsto b$, and let us denote it by $\mcl{B}z$. Then $E$ decomposes as $\mcl{B}z\oplus (\mcl{B}z)^{\perp}$.

  $(ii)\Rightarrow (iii)$: Without loss of generality, we may take $E=\mcl{B}\oplus F$. Then $x\in E$ decomposes as $x=c\oplus y$ with $c\in\mcl{B},~y\in F$. Clearly, $\vphi(b)=\ip{x,bx}=c^*bc+\ip{y,by}$, and we can take $\psi(b)=\ip{y,by}$ for all $b\in\mcl{B}$. Note that $\mcl{B}$ is the two sided (strongly closed) ideal generated by $c$.

  $(iii)\Rightarrow (i)$: Note that the CP-map $b\mapsto c^*bc$ is dominated by the CP-map $\vphi$, and hence there exists a vector $z\in E$ (see \cite{BhatSkeide,paschke73}) such that $c^*bc=\ip{z,bz}$ for all $b\in\mcl{B}$. Note that, for elements $a,a',b,d,d'\in\mcl{B}$, $(acd)^*b(a'cd')=d^*(c^*a^*ba'c)d=d^*\ip{z,a^*ba'z}d'=\ip{azd,ba'zd'}$. It follows that for any element $d$ in the (strongly closed) ideal generated by $c$, there exists an element $z_d\in E$ such that $d^*bd=\ip{z_d,bz_d}$. Taking $d=1$, we have an element $w\in E$ such that $b=\ip{w,bw}$ for all $b\in\mcl{B}$. Observe that $w$ is a unit vector. Direct computation yields $\ip{bw-wb,bw-wb}=0$, hence $w$ is in the center $C_{\mcl{B}}(E)$.
 \eproof


 \bthm \label{rigidity}
  Let $\vphi:\mcl{B}\to\mcl{B}$ be a CP-map such that $\beta(id.,\vphi)<1.$
  Let $(E,x)$ be a GNS-construction for $\vphi$. Then  $E\cong\mcl{B}\oplus F$ for some von Neumann $\mcl{B}$-$\mcl{B}$-module $F$.
 \ethm

 \bproof
  Without loss of generality, assume that $(E,x)$ is the minimal GNS-construction for $\vphi$. Let $\varepsilon>0$  be such that $\beta(id.,\vphi)+\varepsilon< 1$. Since the identity map has $(\mcl{B},1)$ as its GNS-construction,  from Theorem \ref{betainfmod}, there exists $z_1=1\oplus 0, z_2=c\oplus y$ in $\mcl{B}\oplus E$ such that $\norm{z_1-z_2}\leq\beta(id,\vphi)+\varepsilon <1$ and $\vphi(b)=\ip{z_2,bz_2}=c^*bc +\ip{y,by}$.  Further, as
    $\norm{1-c}\leq\norm{z_1-z_2}< 1$
  we note that  $c$ is invertible. Therefore the ideal generated by $c$ is whole of $\mcl{B}$. Now the result follows from the previous Proposition.
 \eproof

 \backn
  The authors wish to thank G. Ramesh and K. B. Sinha for several useful discussions on Bures distance.
 \eackn

  \vspace{.5cm}
  \noindent\sc{Statistics and Mathematics Unit, Indian Statistical Institute, \\R. V. College Post,
  Bangalore 560059, India.}
  \begin{verbatim}E-mail: bhat@isibang.ac.in, sumesh@isibang.ac.in\end{verbatim}

  \begin{center}
   \line(1,0){250}
  \end{center}
\end{document}